# Variants of the Finite Element Method for the Parabolic Heat Equation: Comparative Numerical Study

Ahmed A. Hamada[1,2], Mahmoud Ayyad[1;3], and Amr Guaily[3;4]

[1]Civil, Environmental and Ocean Engineering, Schaefer School of Engineering and Science, Stevens Institute of Technology, Hoboken NJ 07030, USA,

[2]Aerospace and Aeronautical Engineering Department, Faculty of Engineering, Cairo University, Giza 12613, Egypt,
`AHamada@stevens.edu, Ahmed-Hamada03@hotmail.com,`

[3]Engineering Mathematics and Physics Department, Faculty of Engineering, Cairo University, Giza 12613, Egypt,

[4]Smart Engineering Systems Research Center (SESC), School of Engineering and Applied Sciences, Nile University, Shaikh Zayed City,12588, Egypt.

**Abstract.** Different variants of the method of weighted residual finite element method are used to get a solution for the parabolic heat equation, which is considered to be the model equation for the steady state Navier-Stokes equations. Results show that the Collocation and the Least-Squares variants are more suitable for first order systems. Results also show that the Galerkin/Least-Squares method is more diffusive than other methods, and hence gives stable solutions for a wide range of Péclet numbers.

## 1    Introduction

The behavior of different variants of the finite element method is examined using the second order Differential Equation (DE) represented by the parabolic heat equation [1]. The most important features of Finite Element Method (FEM) [2] can be summarized that FEM can approximate complicated geometrical boundaries easily, also it accounts for boundary conditions in an easy, straightforward manner and FEM is modular where the problem may be changed to a more sophisticated problem without changing the program structure. Thus, the used method in this study is the Finite Element Method (FEM).

The weighted residual FEM [3, 4, 5, 6] succeeded to solve boundary value problems for parabolic Differential Equations (DE). The objective of this paper is to study the behaviour of different Finite Element techniques with respect to the second order ordinary DE from two points of view; firstly, the inherent numerical viscosity and the stability of the method with convection dominated flows.

The used FEM and their acronyms that are used throughout the paper are:

- SG: Standard Galerkin method [7],
- C: Collocation method [8, 9, 10],
- LS: Least-Squares method [11, 12, 13, 14, 15, 16],
- GLS: Galerkin/Least-Squares method [17, 18, 19],
- CG: Collocation/Galerkin method [20, 21, 22],
- CLS: Collocation/Least-Squares method [23, 24, 25],
- CGLS: Collocation/Galerkin/Least-Squares method.

The SG and the C methods both depend on the weighted residual method [6]. In the SG method, the weight function is chosen to be a basis function. While, the weight function of the C method is the Dirac delta. C method is used at least since the thirties [10] by Lanczos [8]. The C method is well-known in chemical engineering [3, 6]. The numerical results from the C method is found to be sensitive to the chosen collocation points. The roots of the Legendre polynomials are found to be the optimum collocation points [26, 27, 28].

The LS method is considered as an alternative to the SG formulation. The weight function in the LS is the differential operator of the DE applied on the basis function. The LS FEM is receiving high attention in different communities [29, 30, 31, 32, 33]. In 1988, Hughes and Shakib [17, 18, 19] has presented the GLS method to solve the advective diffusive systems. The GLS is based on adding the differential operator, introduced in the LS, to the SG method to improve the stability of the advective diffusive equation [34].

The paper is organized as follows. The differential equation is described in Section 2. Then, the exact and the numerical solutions, using the finite element method, of the differential equation are discussed in section 3. Finally, section 4 presents the discussion of the results and the conclusion.

## 2    Problem Statement

Common example of one dimensional (1D) second order differential equations is the parabolic heat equation. The DE is solved numerically using the weighted residual Finite Element Method (FEM). The non-dimensional DE is expressed in the Cartesian coordinate as,

$$T_{xx}(x) = P_e T_x(x) \qquad (1)$$

where, $T(x)$ is the Temperature at distance $x$, the subscript $x$ denotes the derivative with respect to $x$, $P_e$ is the Péclet number ($P_e = \frac{Lu}{\alpha}$; $L$ is the characteristic length of the rod, $u$ is the local flow velocity, and $\alpha$ is the thermal diffusivity). Dirichlet boundary conditions are used and are discussed later in section 3.2.2.

# 3 Solving the Differential Equation

In this section, the exact and numerical solution of the parabolic heat equations 1 is discussed.

## 3.1 The Exact Solution

The Exact solution of the parabolic heat equation 1 is found for the used boundary conditions; $T(1) = 0$ and $T(2) = 1$, using the principle of superposition,

$$T(x) = \frac{e^{P_e(x-1)} - 1}{e^{P_e} - 1} \tag{2}$$

## 3.2 The Finite Element Solution

### 3.2.1 The Weak Form

The weak forms of the parabolic heat equation 1 is,

$$\int_{-1}^{1} [(P_e T_x(x) - T_{xx}(x))w_i(x)] dx = 0 \tag{3}$$

where, $w_i(.)$ are the weight functions. Integration by parts is applied to the higher order derivative term in equation 3. The generated term, product of functions, is neglected because of using Dirichlet boundary condition which yields to,

$$\int_{-1}^{1} \left[ P_e T_x(x)w_i(x) + T_x(x)w_{i_x}(x) \right] dx = 0 \tag{4}$$

The linear shape functions, $N_j(.)$, are used to approximate the dependent quantity $T(x) = N_j(x)T_j$ between nodes, where $T_j$ are the nodal values. Then, the weak form of equation 4 can be rewritten in the indicial notation as,

$$\int_{-1}^{1} \left[ P_e N_{j_x}(x)w_i(x) + N_{j_x}(x)w_{i_x}(x) \right] dx T_j = 0 \tag{5}$$

### 3.2.2 The Weight Function

The used weight function depends on the used FE technique. As mentioned in the introduction, section 1, different techniques are used to compare the performance of the FE techniques. The weight function in the case of the SG is the shape functions themselves, equation 6a, in the C method is the Dirac delta function, equation 6b, and in the LS method is a differential operator of the DE, equation 6c. The GLS, CG, CLS, and the CGLS are combinations of the main three methods, SG, C, and LS, as shown in equations 6d-6g, respectively.

$$w_i(x) = N_i(x) \tag{6a}$$

$$w_i(x) = \delta_i(x) \tag{6b}$$

$$w_i(x) = L(N_i(x)) \tag{6c}$$

$$w_i(x) = N_i(x) + \tau L(N_i(x)) \tag{6d}$$

$$w_i(x) = \delta_i(x) + N_i(x) \tag{6e}$$

$$w_i(x) = \delta_i(x) + \tau L(N_i(x)) \tag{6f}$$

$$w_i(x) = \delta_i(x) + N_i(x) + \tau L(N_i(x)) \tag{6g}$$

where, $\delta_i(x)$ is the Dirac delta, $\tau$ is the stabilization parameter and $L(N_i(x))$ is the differential operator. $\tau$ is the stabilization parameter given by [18],

$$\tau = \left( \left( \frac{2P_e}{l_e} \right)^2 + \left( \frac{4}{l_e^2} \right)^2 \right)^{-0.5} \tag{7}$$

where, $l_e$ is the length of the element in a uniform grid. $L(N_i(x))$ is the differential operator given by,

$$L(N_i(x)) = P_e N_{i_x}(x) - N_{i_{xx}}(x) \tag{8}$$

The integration is calculated using the Gaussian Quadrature points [35]. The first derivative of the Dirac delta function, using its properties, is,

$$\delta_x(x) = \frac{-\delta(x)}{x} \tag{9}$$

## 4    Results and Discussion

### 4.1    The Exact Solution

The exact solution is a straight line at a zero Péclet number (a pure diffusion). As the Péclet number increases, the exact solution of the parabolic heat equation deviates from the straight line, as a consequence of the equation being convection dominated, as shown in figure 1. As a result, a boundary layer starts to form and the higher the Péclet number the thinner the boundary later.

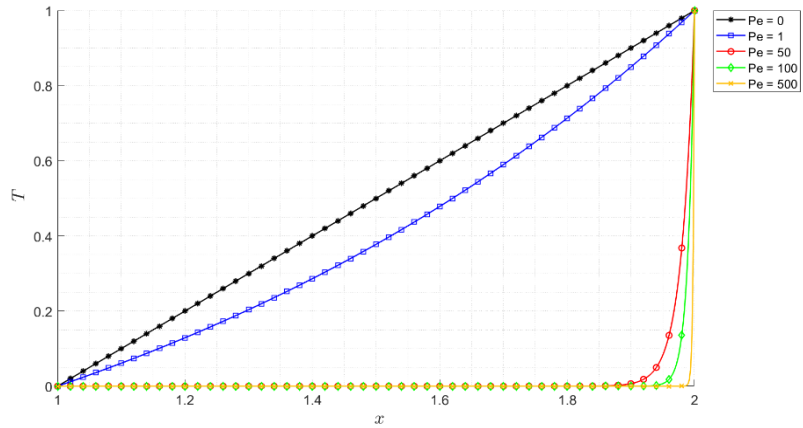

**Fig. 1.** Exact solutions of the parabolic heat equation at different Péclet numbers.

### 4.2 The Finite Element Solution

The numerical solution of the parabolic heat equation 5 is conducted for different number of elements; 25, 50, 75, and 100 at different Péclet Numbers; 1, 50, 100, and 500 by the techniques mentioned previously. The adopted procedure, integration by parts, to deal with the higher derivative term in the DE produces a singular matrix in the C method with very small condition number. Thus, the solution from the C method ceased to exist. Figures 2-7 show the effect of changing the Péclet numbers for each method. Figures 2 and 3 show that the SG and CG methods become unstable as the Péclet number increases because the effect of the dispersive term increases. Considering linear shape functions for the second order DE leads to the conclusion that the LS method to become a pure diffusive method whatever the value of Péclet number as shown in figure 4. The reasoning behind this behavior for the LS method is that the weight function is the derivative of the residual of the equitation which contains second derivative of the linear shape functions and consequently vanishes. The remedy is to use shape functions of the same order as the highest derivative appearing in the DE.

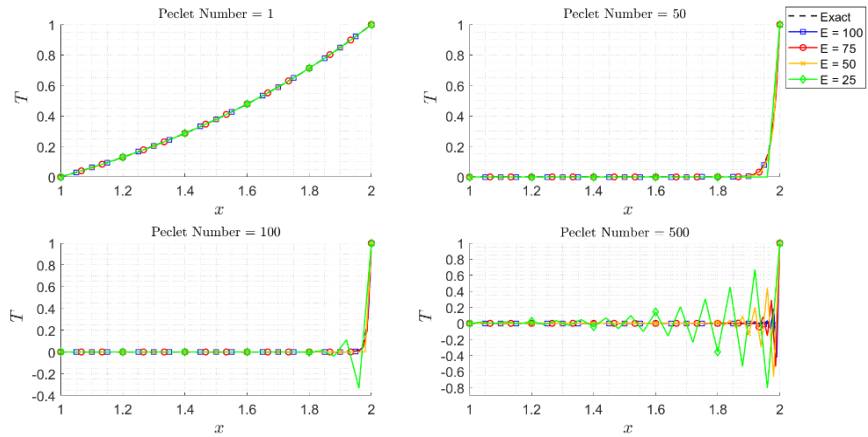

**Fig. 2.** Exact and numerical solutions using the SG FEM for the parabolic heat equation at different number of elements and different Péclet numbers.

The CLS becomes more stable as the Péclet number increases, as shown in figure 5, because the effect of the diffusive term from the LS increases. Figures 6 and 7 show that the GLS method is always stable whatever the value of the Péclet number, while the CGLS method becomes unstable at higher Péclet number because the effect of the convection term in the SG part is higher than the effect of the diffusivity in the LS part. Moreover, figure 8 shows that the GLS is the most accurate diffusive method. Hence, the results proved that GLS is the most recommended method to solve second order DE, such as Navier-Stokes equations, with high accuracy [36, 37]. However, solving the first order DE, such as Euler equation or Navier-Stokes equations formulated as first order system [38], with LS and C gives a solution with higher accuracy than other methods [10].

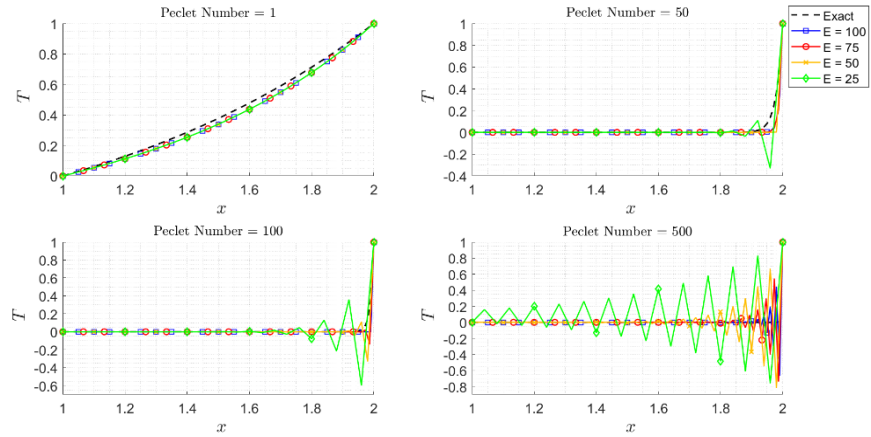

**Fig. 3.** Exact and numerical solutions using the CG FEM for the parabolic heat equation at different number of elements and different Péclet numbers.

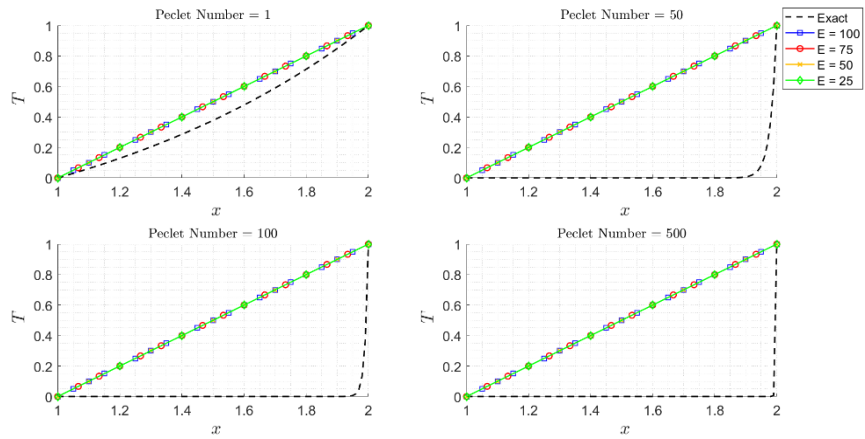

**Fig. 4.** Exact and numerical solutions using the LS FEM for the parabolic heat equation at different number of elements and different Péclet numbers.

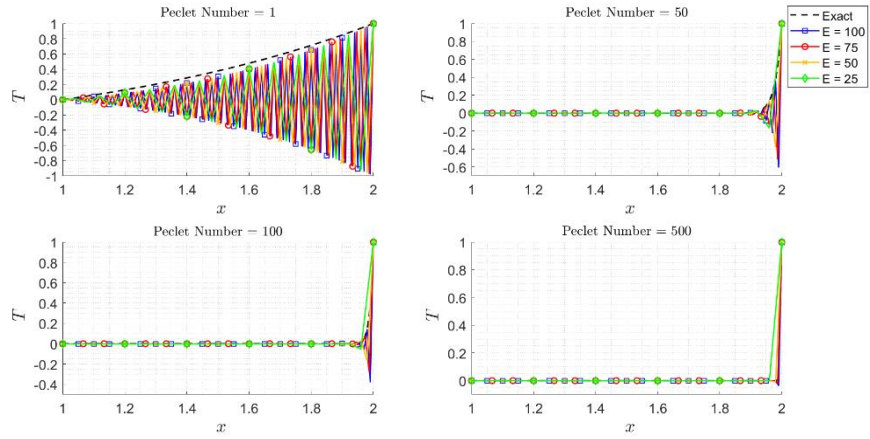

**Fig. 5.** Exact and numerical solutions using the CLS FEM for the parabolic heat equation at different number of elements and different Péclet numbers.

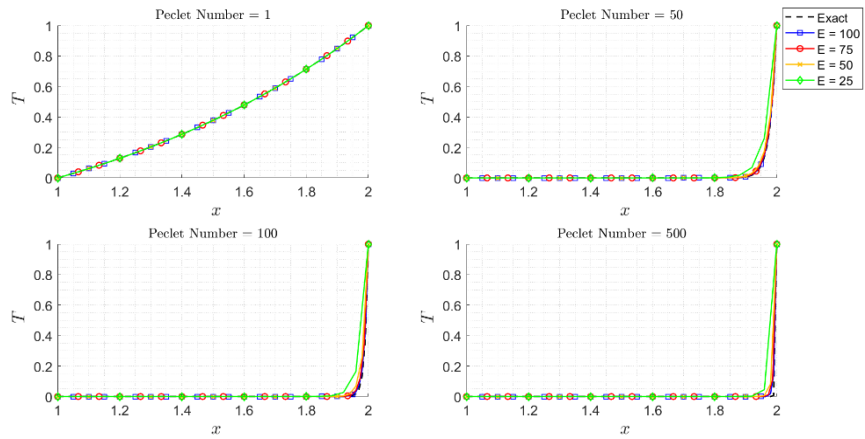

**Fig. 6.** Exact and numerical solutions using the GLS FEM for the parabolic heat equation at different number of elements and different Péclet numbers.

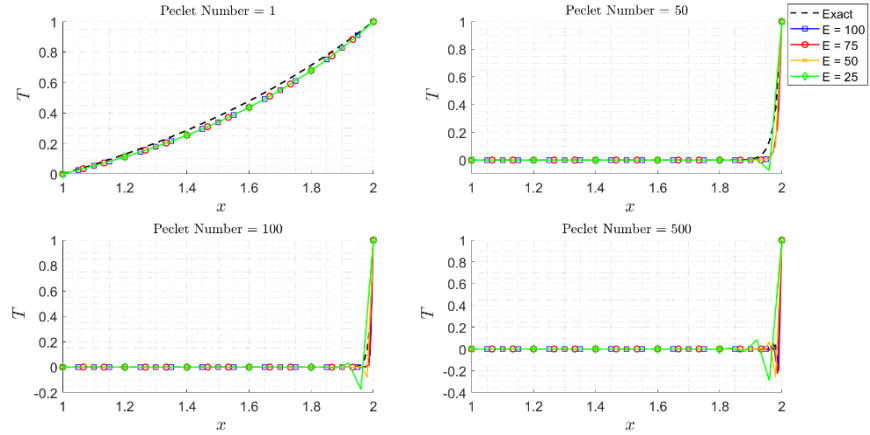

**Fig. 7.** Exact and numerical solutions using the CGLS FEM for the parabolic heat equation at different number of elements and different Péclet numbers.

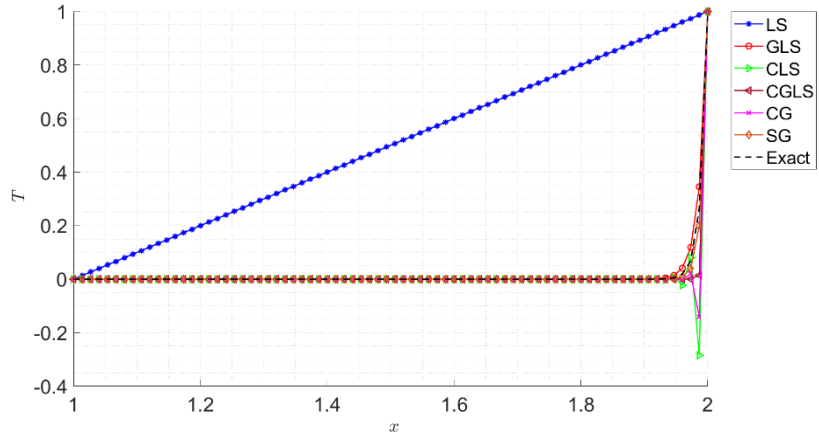

**Fig. 8.** The exact solution at Péclet number 100 and the numerical solutions using 75 elements.

Figures 9 and 10 show the accuracy and the convergence of the FE methods using the absolute relative error,

$$\text{Absolute Relative Error} = \left| \frac{\text{Exact Solution} - \text{Numerical Solution}}{\text{Exact Solution}} \right| \qquad (10)$$

The GLS method gives high accuracy solution regardless the value of the Péclet number and converges by increasing the number of elements. The solution from SG method shows promising solution at small values of Péclet number and the errors decrease with the increase of number of elements. However, the solution from the SG methods diverges when the value of Péclet number is increased. Moreover, the methods CLS, CG, and CGLS have nearly the same error because, as discussed previously, the dispersive term effects are higher than the diffusion term effects. Finally, these figures emphasize that the LS method fails to solve the second order DE with linear shape functions.

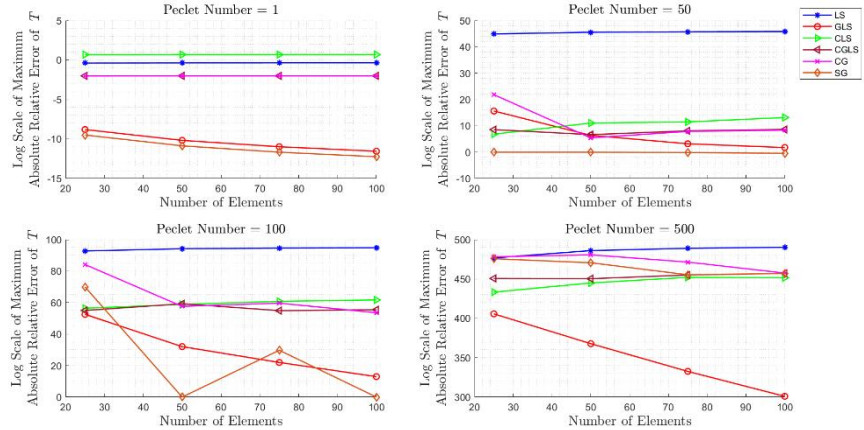

**Fig. 9.** The Logarithmic scale of maximum absolute relative error of the numerical solutions at different number of elements and different Péclet numbers.

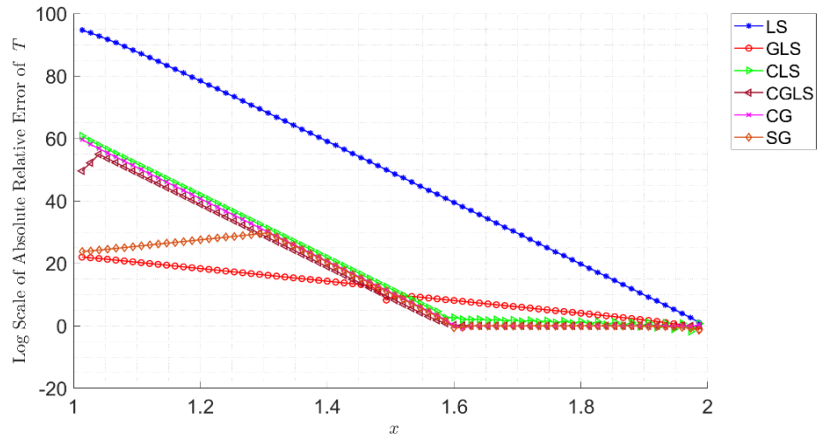

**Fig. 10.** The Logarithmic scale of absolute relative error of the numerical solutions at Péclet number 100 using 75 elements.

## 5 Conclusion

The parabolic heat equation is used to study the behavior of different versions of the weighted residuals method. The Péclet number is varied from 1, which is considered to be diffusion dominated case, till a very high value of about 500, reflecting the case of convection dominated case. The GLS method turned out to be the most stable technique for the case at hand, while the CLS method is believed to be suitable for very high Péclet number case or in general for the inviscid flows which correspond to very high Reynolds number flows.